\documentclass[11pt, numbers=enddot, twoside, headsepline, abstract=true]{scrartcl}
\usepackage[margin=1in, a4paper]{geometry}

\usepackage{etoolbox}
\setkomafont{title}{\sffamily\bfseries}
\makeatletter
\patchcmd{\@maketitle}{\huge}{\LARGE}{}{}
\makeatother

\usepackage{authblk}

\makeatletter
\let\oldAB@maketitle\AB@maketitle
\renewcommand{\AB@maketitle}{%
  \oldAB@maketitle
  \vspace{-4.5em} 
}
\makeatother

\usepackage{indentfirst}

\usepackage[numbers,sort&compress]{natbib}
\usepackage[hypertexnames=false,colorlinks,linkcolor={blue},citecolor={blue},urlcolor={blue}]{hyperref}

\usepackage{amssymb}
\usepackage{amsmath}
\usepackage{stmaryrd}
\usepackage{empheq}
\usepackage{amsthm}
\usepackage[mono=false]{libertine}
\usepackage{mathrsfs}
\usepackage[libertine, 
nosymbolsc, nonewtxmathopt, subscriptcorrection, vvarbb]{newtxmath}
\usepackage{euscript}

\usepackage{graphicx}
\usepackage{epsfig}
\usepackage{float}
\usepackage{caption}
\usepackage{subcaption}
\usepackage{tabularray}

\usepackage{enumitem}
\setlist{
	leftmargin=*, 
	topsep = 4pt, 
	itemsep = 1pt, 
  wide = 0pt
}

\numberwithin{equation}{section}

\theoremstyle{definition}

\theoremstyle{remark}

\allowdisplaybreaks

\newcommand{\reDeclareMathOperator}[2]{\let#1\undefined \DeclareMathOperator{#1}{#2}}
\newcommand{\reDeclareMathOperatorL}[2]{\let#1\undefined \DeclareMathOperator*{#1}{#2}}

\newcommand{\eqfinp}{\,.}
\newcommand{\eqfinv}{\,,}

\usepackage{macros}

\usepackage{tikzit}

\tikzstyle{title fill}=[-, fill={rgb,255: red,223; green,223; blue,223}]
\tikzstyle{superblock}=[-, fill=none, dashed]
\tikzstyle{arrow}=[->]

\title{Optimal class scheduling with minimum distributed vacant seats using binary linear programming}

\author[1]{Ohmchana Klinsod}
\author[1,$\star$]{Parin Chaipunya}
\affil[1]{Department of Mathematics, Faculty of Science, King Mongkut's University of Technology Thonburi,\protect\\ 126 Pracha Uthit Rd., Bang Mod, Thung Khru, Bangkok 10140, Thailand}
\affil[$\star$]{Corresponding author.\ Email: parin.cha@kmutt.ac.th}
\date{}

\def\runningauthor{O. Klinsod and P. Chaipunya}
\def\runningtitle{Optimal class scheduling with minimum distributed vacant seats...}

\usepackage{scrlayer-scrpage}
\clearpairofpagestyles
\setkomafont{pagehead}{\normalfont\small}
\ofoot{}
\lehead{\runningauthor}
\lohead{\runningtitle}
\rehead{\pagemark}
\rohead{\pagemark}

\begin{document}
\maketitle 

\begin{abstract}
  In this paper, we propose a mathematical model for class scheduling using binary linear programming with the aims to automate the tedious task and to minimize the leftover seats distributed over all the classrooms.
  A case study at the Department of Mathematics, King Mongkut's University of Technology Thonburi (KMUTT), Thailand, is also presented.
  The result from our model is shown to have improved significantly over the current manual scheduling method.

	
	
	
\end{abstract}


\normalsize

{
  \tableofcontents
}

\section{Introduction}

Class scheduling in any academic institution has been a notoriously complex combinatoric and time-consuming task.
Most of the scheduling jobs are carried out manually, yielding a feasible but suboptimal schedule which may lack the efficiency in terms of resource utilization.
The challenge of automated scheduling lies in the simultaneos allocation of deeply connected resrouces including instructing facility, student cohorts, temporal slots, and physical infrastructure requirements such as computers or chemical substances for laboratories.
Due to resource limitation as well as other physical and operational constraints, developing an optimal schedule manually can be far from being optimal.
To mitigate these inefficiencies, rigorous mathematical frameworks have been developed in order to construct an optimal timetable that maximizes resource or convenience criteria.
These decisive factors may be measured by, {\itshape e.g.}, the spatial mismatches, resting times, or travel times, under practical physical and operational constraints.

This scheduling problem lies within the broader framework of operations research, discrete optimization and scheduling theory.
The general modeling principle using optimization-based techniques has been discussed extensively in standard texts such as \citet{zbMATH04045470}, while the dedicated works scheduling theory may be found in, {\itshape e.g.}, \citet{zbMATH07555488}.
More specialized studies of class timetabling has also been widely studied as a constrained optimization problem involving the assignment of courses, instructors, student groups, rooms, and time periods subject to resource, availability, and operational restrictions.
In the work of \citet{zbMATH01746207}, the authors highlighted the use of heuristics and evolutionary algorithm, decomposition ideas for large real-world timetabling instances, multicriteria decision-making approaches, and case-based reasoning. 
This is particularly relevant to the present work because it shows that practical timetabling is more than just finding a feasibility schedule, but also a problem of balancing several competing quality criteria.
Similarly, \citet{zbMATH05258612} surveys metaheuristic techniques for university timetabling, including graph-based heuristics, local search, tabu search, simulated annealing, evolutionary methods, and hybrid approaches. 
These methods are especially important for large instances, where exact optimization models may become computationally demanding.
For more recent advances of class timetabling, the surveys by \citet{zbMATH07709092}, \citet{9499056} and works of \citet{zbMATH06736101} and \citet{zbMATH05634685} provide plenty of useful information especially on the solution approaches.

Apart from the development of solution techniques, heuristics, and compuational tricks for class timetabling problem, there are also majority of researches that are dedicated to incorporating and implementing complicate practical features into the scheduling model.
These involves the formulation of new objective functions and crafting the physical and operational needs into equality and inequality constraints.
In the work of \citet{zbMATH06635852}, the authors allow the teachers' preferences to influence the resulting schedule to promote satisfactory of the teaching faculties.
Later \citet{zbMATH08018028} studied a similar problem but also consider the room stability attribute of a timetable.
This means that the model tries to assign a group of students to the same room as much as possible, minimizing the room rotation from the point of view of students.
On the other hand, \citet{zbMATH08066550} considered classes that are conducted on-campus, online, and hybrid of the two modes.
This work used a multi-objective binary programming framework to find an optimal solution with the right balances between resource utilization, student movement, and class modes.
Other metrics like fairness was considered by \citet{zbMATH06581562}.

Although class and course timetabling have been extensively studied in the literature, practical implementations often require institution-specific constraints that go beyond what most of the benchmark and existing frameworks have provided.
In this work, we therefore formulate a model that explicitly treats several important features that is crucial in practice.
First, we impose type-matching constraints that link course metadata directly to physical infrastructure. 
Lecture-based subjects must be assigned to lecture rooms, while laboratory-based subjects must be assigned to laboratory rooms of the same nature. 
This promotes the appropriate use of resources and helps avoiding the need for post-optimization manual adjustments.

We also define synchronized parallel-course constraints that manage multi-section courses requiring simultaneous delivery.
This logic is useful in coordinated-teaching activities, or common assessment sessions.
By imposing this directly into the problem's constraints, the resulting timetable respects both physical classroom and temporal synchronization requirements.

Further more, we design the objective function to penalize capacity mismatch that is sensitive to the distribution.
For instance, the purely linear metric that only measure the total empty seats does not distinguish between (a) having~$10$ vacant seats in a single room, and (b) having two rooms with~$5$ vacant seats each.
By using the squared empty seat penalization, the model places more severe cost on stronger mismatches and hence encouraging a more balanced solution.
Even the squared penalty is used, the objective remains linear in its decision variables.
Lastly, to combat the combinatorial explosion and NP-hard complexity typical of large institutional datasets, we introduce a preprocessing phase that partitions individual students into exact cohort equivalence classes based on shared enrollment requirements, significantly reducing the runtime and dimensional scale of the solver.

To support our modeling design, we also include two sample implementation of our models.
First, we begin with a toy example which serves as an illustrative demonstration for the readers to absorb how our model is used.
Then we present a real case study, implementing the model with the real data from the Department of Mathematics, King Mongkut's University of Technology Thonburi (KMUTT), for its class timetabling of Semester 1/2025.
In both implementations, the model is solved exactly using the solver and preconditioning of {\ttfamily Gurobi Optimizer} version 13 (with academic license) via the use of {\ttfamily Python} interface {\ttfamily gurobipy}.

\section{Methodology}

We start the section with~\S\ref{subsec: problem characterization} to describe the problem and how the data for our main case study is processed.
The problem formulation as well as all the requirements are discussed in~\S\ref{subsec: model}.

\subsection{Problem characterization and data pre-processing}
\label{subsec: problem characterization}

The class timetabling problem examined in this study is based on the operational and physical restrictions used at the Department of Mathematics at King Mongkut's University of Technology Thonburi (KMUTT).
Its dataset spanning the second semester of the 2025 academic year is also used for our main case study in~\S\ref{sec: case study}.
The underlying scheduling environment consists of five primary discrete sets (see also Table~\ref{tab:notations}), namely Courses, Temporal slots, Classrooms, Instructors, and Student cohorts.

To transform the problem into a clean mathematical layout, we develop an automation script that extract and parse the raw data into a structured one using the {\ttfamily Python}'s {\ttfamily pandas} library.
This data pipeline handles initial formatting, flags missing fields, and draw the required parameters including the classroom capacities, instructor availability, and other requirements.

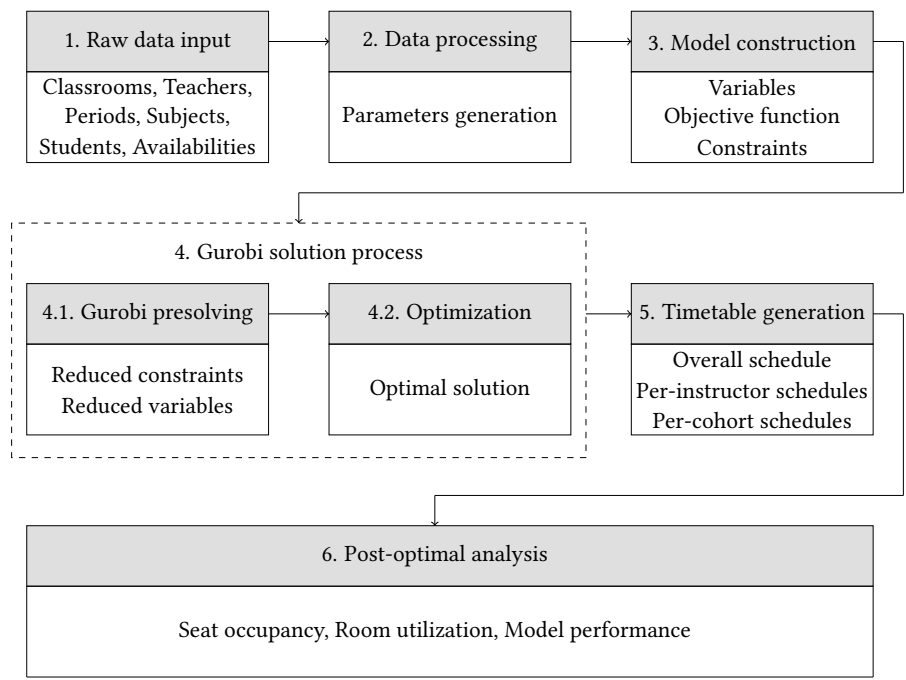
\begin{figure}[H]
  \centering
  \scalebox{0.8}{\begin{tikzpicture}
	\begin{pgfonlayer}{nodelayer}
		\node [style=none] (0) at (-8, 2) {};
		\node [style=none] (1) at (0, 2) {};
		\node [style=none] (2) at (-8, 0) {};
		\node [style=none] (3) at (0, 0) {};
		\node [style=none] (4) at (-4, 1) {1. Raw data input};
		\node [style=none] (5) at (2, 2) {};
		\node [style=none] (6) at (10, 2) {};
		\node [style=none] (7) at (2, 0) {};
		\node [style=none] (8) at (10, 0) {};
		\node [style=none] (9) at (6, 1) {2. Data processing};
		\node [style=none] (10) at (12, 2) {};
		\node [style=none] (11) at (20, 2) {};
		\node [style=none] (12) at (12, 0) {};
		\node [style=none] (13) at (20, 0) {};
		\node [style=none] (14) at (16, 1) {3. Model construction};
		\node [style=none] (15) at (-8, -7) {};
		\node [style=none] (16) at (0, -7) {};
		\node [style=none] (17) at (-8, -9) {};
		\node [style=none] (18) at (0, -9) {};
		\node [style=none] (19) at (-4, -8) {4.1. Gurobi presolving};
		\node [style=none] (20) at (-4, -0.5) {Classrooms, Teachers,};
		\node [style=none] (21) at (-4, -1.5) {Periods, Subjects,};
		\node [style=none] (22) at (-4, -2.5) {Students, Availabilities};
		\node [style=none] (23) at (-8, 0) {};
		\node [style=none] (24) at (-8, -3) {};
		\node [style=none] (25) at (0, 0) {};
		\node [style=none] (26) at (0, -3) {};
		\node [style=none] (27) at (6, -1.5) {Parameters generation};
		\node [style=none] (30) at (2, 0) {};
		\node [style=none] (31) at (2, -3) {};
		\node [style=none] (32) at (10, 0) {};
		\node [style=none] (33) at (10, -3) {};
		\node [style=none] (34) at (16, -0.5) {Variables};
		\node [style=none] (35) at (16, -1.5) {Objective function};
		\node [style=none] (36) at (16, -2.5) {Constraints};
		\node [style=none] (37) at (12, 0) {};
		\node [style=none] (38) at (12, -3) {};
		\node [style=none] (39) at (20, 0) {};
		\node [style=none] (40) at (20, -3) {};
		\node [style=none] (41) at (2, -7) {};
		\node [style=none] (42) at (10, -7) {};
		\node [style=none] (43) at (2, -9) {};
		\node [style=none] (44) at (10, -9) {};
		\node [style=none] (45) at (6, -8) {4.2. Optimization};
		\node [style=none] (46) at (-4, -10) {Reduced constraints};
		\node [style=none] (47) at (-4, -11) {Reduced variables};
		\node [style=none] (49) at (-8, -9) {};
		\node [style=none] (50) at (-8, -12) {};
		\node [style=none] (51) at (0, -9) {};
		\node [style=none] (52) at (0, -12) {};
		\node [style=none] (53) at (6, -10.5) {Optimal solution};
		\node [style=none] (54) at (2, -9) {};
		\node [style=none] (55) at (2, -12) {};
		\node [style=none] (56) at (10, -9) {};
		\node [style=none] (57) at (10, -12) {};
		\node [style=none] (58) at (12, -7) {};
		\node [style=none] (59) at (20, -7) {};
		\node [style=none] (60) at (12, -9) {};
		\node [style=none] (61) at (20, -9) {};
		\node [style=none] (62) at (16, -8) {5. Timetable generation};
		\node [style=none] (63) at (16, -9.5) {Overall schedule};
		\node [style=none] (64) at (16, -10.5) {Per-instructor schedules};
		\node [style=none] (65) at (16, -11.5) {Per-cohort schedules};
		\node [style=none] (66) at (12, -9) {};
		\node [style=none] (67) at (12, -12) {};
		\node [style=none] (68) at (20, -9) {};
		\node [style=none] (69) at (20, -12) {};
		\node [style=none] (70) at (-8, -15) {};
		\node [style=none] (71) at (20, -15) {};
		\node [style=none] (72) at (-8, -17) {};
		\node [style=none] (73) at (20, -17) {};
		\node [style=none] (74) at (5.5, -16) {6. Post-optimal analysis};
		\node [style=none] (75) at (5.5, -18.5) {Seat occupancy, Room utilization, Model performance};
		\node [style=none] (76) at (-8, -17) {};
		\node [style=none] (77) at (-8, -20) {};
		\node [style=none] (78) at (20, -17) {};
		\node [style=none] (79) at (20, -20) {};
		\node [style=none] (80) at (-8.5, -5) {};
		\node [style=none] (81) at (-8.5, -12.75) {};
		\node [style=none] (82) at (10.5, -5) {};
		\node [style=none] (83) at (10.5, -12.75) {};
		\node [style=none] (84) at (1, -6) {4. Gurobi solution process};
		\node [style=none] (85) at (0, 1) {};
		\node [style=none] (86) at (2, 1) {};
		\node [style=none] (87) at (10, 1) {};
		\node [style=none] (88) at (12, 1) {};
		\node [style=none] (89) at (1, -5) {};
		\node [style=none] (90) at (1, -4) {};
		\node [style=none] (91) at (21, -4) {};
		\node [style=none] (92) at (21, 1) {};
		\node [style=none] (93) at (20, 1) {};
		\node [style=none] (94) at (10.5, -8) {};
		\node [style=none] (95) at (12, -8) {};
		\node [style=none] (96) at (20, -8) {};
		\node [style=none] (97) at (21, -8) {};
		\node [style=none] (98) at (21, -14) {};
		\node [style=none] (99) at (5.5, -15) {};
		\node [style=none] (100) at (5.5, -14) {};
		\node [style=none] (101) at (0, -8) {};
		\node [style=none] (102) at (2, -8) {};
	\end{pgfonlayer}
	\begin{pgfonlayer}{edgelayer}
		\draw [style=superblock] (82.center)
			 to (80.center)
			 to (81.center)
			 to (83.center)
			 to cycle;
		\draw [style=title fill] (0.center)
			 to (1.center)
			 to (3.center)
			 to (2.center)
			 to cycle;
		\draw [style=title fill] (5.center)
			 to (6.center)
			 to (8.center)
			 to (7.center)
			 to cycle;
		\draw [style=title fill] (10.center)
			 to (11.center)
			 to (13.center)
			 to (12.center)
			 to cycle;
		\draw [style=title fill] (15.center)
			 to (16.center)
			 to (18.center)
			 to (17.center)
			 to cycle;
		\draw (25.center)
			 to (23.center)
			 to (24.center)
			 to (26.center)
			 to cycle;
		\draw (32.center)
			 to (30.center)
			 to (31.center)
			 to (33.center)
			 to cycle;
		\draw (39.center)
			 to (37.center)
			 to (38.center)
			 to (40.center)
			 to cycle;
		\draw [style=title fill] (41.center)
			 to (42.center)
			 to (44.center)
			 to (43.center)
			 to cycle;
		\draw (51.center)
			 to (49.center)
			 to (50.center)
			 to (52.center)
			 to cycle;
		\draw (56.center)
			 to (54.center)
			 to (55.center)
			 to (57.center)
			 to cycle;
		\draw [style=title fill] (60.center)
			 to (58.center)
			 to (59.center)
			 to (61.center)
			 to cycle;
		\draw (68.center)
			 to (66.center)
			 to (67.center)
			 to (69.center)
			 to cycle;
		\draw [style=title fill] (70.center)
			 to (71.center)
			 to (73.center)
			 to (72.center)
			 to cycle;
		\draw (78.center) to (76.center);
		\draw (76.center) to (77.center);
		\draw (77.center) to (79.center);
		\draw (79.center) to (78.center);
		\draw [style=arrow] (85.center) to (86.center);
		\draw [style=arrow] (87.center) to (88.center);
		\draw (93.center) to (92.center);
		\draw (92.center) to (91.center);
		\draw (91.center) to (90.center);
		\draw [style=arrow] (90.center) to (89.center);
		\draw [style=arrow] (94.center) to (95.center);
		\draw (96.center) to (97.center);
		\draw (97.center) to (98.center);
		\draw (98.center) to (100.center);
		\draw [style=arrow] (100.center) to (99.center);
		\draw [style=arrow] (101.center) to (102.center);
	\end{pgfonlayer}
\end{tikzpicture}}
  \caption{The workflow and data pipeline of our modeling process.}
  \label{fig: process}
\end{figure}

After the required data is extracted from the raw sheet, it is fed into the model as a set of parameters. 
The model is built accordingly to the formulation in~\S\ref{subsec: model} and solved using {\ttfamily Gurobi Optimizer}.
The optimal solution is then parsed into a {\ttfamily Python} script to generate an overall schedule, per-instructor schedules, and per-cohort schedules.
This script uses the {\ttfamily pandas} package with the help of {\ttfamily openpyxl} engine.
The full workflow is visualized as in Figure~\ref{fig: process}.

One of the challenges in solving exactly an integer program is the exponential-order complexity of the integer variables.
To overcome this bottlenect, this study introduces a processing phase that maps individual student enrollment data to equivalence classes, called \emph{cohorts}.
For instance, let~$G_{\text{bip}} = (U, I , E)$ denote a bipartite graph where~$U$ represents the vertex set of individual students,~$I$ represents the vertex set of available courses, and~$E \subset U \times I$ the set of edge in which~$(u, i) \in E$ if and only if the student~$u$ enrolls in the course~$i$.
This graph~$G_{\text{bip}}$ is completely idenfied by its adjacency matrix~$M(G_{\text{bip}}) = [M_{ui}] \in \R^{\abs{U} \times \abs{I}}$, where
\begin{align*}
  M_{ui} =
  \begin{cases}
    1 &\text{if~$(u,i) \in E$} \eqfinv \\
    0 &\text{otherwise} \eqfinp
  \end{cases}
\end{align*}
Then we define an equivalence relation~$\sim$ on~$U$ by
\begin{align*}
  u_{1} \sim u_{2} \iff M_{u_{1},\bullet} = M_{u_{2},\bullet} 
  \eqfinv
\end{align*}
that is, if their corresponding rows in~$M(G_{\text{bip}})$ equal.
The equivalence classes obtained from~$\sim$ partition~$U$, and these partitions are recognized as \emph{cohorts} of students whose enrollment routines are exactly the same.
The set of all cohorts is denoted by~$G$.
This compression could drasitcally scale down the dimension of the search space.
In particular, applying this graph partitioning technique to our dataset reduces the student dimension from~$489$ down to~$126$ distinct cohorts.

\subsection{Mathematical model formulation}
\label{subsec: model}

This subsection focuses on the mathematical model formulation based on the practice and requirements at the Department of Mathematics, King Mongkut's University of Technology Thonburi (KMUTT).
Before entering the modeling formalism, let us summarize all the involved sets and parameters used in our study in the following Table~\ref{tab:notations}.
\begin{longtblr}[
  caption = {Summary of sets and parameters},
  label = {tab:notations}
  ]{
    colspec = {l Q[l, co=1]},
    row{1} = {font=\bfseries},
    rowhead = 1,
    hline{1,2,Z} = {0.08em},
    rowsep = 5pt,
  }
  Notation & Description \\
  {\bfseries Sets} \\
  $I$ & Set of all scheduled courses, indexed by~$i$. \\
  $J$ & Set of all physical classroom spaces, indexed by~$j$. \\
  $K$ & Set of all discrete feasible temporal slots per week, indexed by~$k$. \\
  $T$ & Set of all faculty instructors, indexed by~$t$. \\
  $G$ & Set of all reduced student cohorts, indexed by~$g$. \\
  $K_{\text{lunch}}$ & Subset of temporal blocks reserved for lunch breaks. \\
  $I_{t} \subset I$ & Subset of courses assigned to instructor~$t$. \\
  $I_{g} \subset I$ & Subset of courses required by student cohort~$g$. \\
  $I^{\text{PAR}}_{t} \subset 2^{I_{t}}$ & A collection of subsets of~$I_{t}$.
  Its element~$I^{\text{par}}_{t} \in I^{\text{PAR}}_{t}$ is a subset of~$I_{t}$ containing courses, taught by~$t \in T$, that have to be parallel to one another.   \\
  $I_{\text{lec}}, I_{\text{lab}} \subset I$ & Disjoint subsets of lecture and laboratory courses, respectively. \\
  $J_{\text{lec}}, J_{\text{lab}} \subset J$ & Disjoint subsets of lecture rooms and laboratories, respectively. \\
  $K_{t}^{\text{unavail}} \subset K$ & Subset of time slots where instructor~$t$ is unavailable. \\
  $K_{g}^{\text{unavail}} \subset K$ & Subset of time slots where student cohort~$g$ is restricted from classes. \\
  $J_{k}^{\text{unavail}} \subset J$ & Subset of physical classrooms pre-allocated to external departments. \\
  {\bfseries Parameters} \\
  $H_{i}$ & Total continuous duration (in hours) required for course~$i$. \\
  $C_{j}$ & Maximum seating or workstation capacity of physical classroom~$j$. \\
  $S_{i}$ & Total enrollment size (number of active students) registered for course~$i$. \\
  $D_{k}$ & Calendar day identifier corresponding to time slot~$k$. \\
  $R$ & Resting hour parameter. \\
\end{longtblr}

\subsubsection{Decision variables and objective function}
To express our allocation decisions, we establish two families of binary decision variables:
\begin{align*}
  x_{ijk} 
  &= 
  \begin{cases} 
    1 & \text{\small if course~$i$ is conducted in classroom~$j$ during time slot~$k$} \eqfinv \\
    0 & \text{\small otherwise} \eqfinv 
  \end{cases}
  \intertext{and}
  y_{ijk} 
  &= 
  \begin{cases} 
    1 & \text{\small if course~$i$ begins its block in classroom~$j$ at the start of time slot~$k$} \eqfinv \\
    0 & \text{\small otherwise} \eqfinp
  \end{cases}
\end{align*}
Although it seems redundant to have both families at first glance, they greatly help expressing the cost function and constraints with better clarity.

The objective function seeks to minimize the total sum of squared empty seats across all the operational configurations:
\begin{align*}
  Z = \sum_{i \in I} \sum_{j \in J} \sum_{k \in K} y_{ijk} (C_{j} - S_{i})^{2} 
  \quad \leftarrow \text{Minimize}
\end{align*}
where~$C_{j}$ represents the total seating capacity of room~$j$, and~$S_{i}$ represents the total enrollment size of the subject~$i$.
Note that this objective remains linear since the squares occur strictly at the coefficients.
Squaring the capacity gap aggressively forces small cohorts away from large lecture halls based on the properties of the \emph{rearrangement inequality}, driving the system toward an optimized spatial distribution.
Without the squares, the objective would have been~$\tilde{Z} = \sum_{i \in I} \sum_{j \in J} \sum_{k \in K} y_{ijk} (C_{j} - S_{i})^{2}$.
This alternative is not a good choice because it does not help distributing the empty seats across different classrooms.
For example,~$\tilde{Z}$ could not distinguish between (a) having one empty seat in each of the~$100$ classrooms, and (b) having~$100$ empty sets in one classroom.

One would observe due to the binary variables, objective function, and further constraints in~\S\ref{subsec: constraints} that finally the model is a binary linear programming (BLP).

\subsubsection{Constraints}\label{subsec: constraints}

In this subsection, we list all the constraints that frame the operational and physical restrictions of class scheduling practice at the Department of Mathematics, KMUTT.
\begin{itemize}[label=$\circ$, leftmargin=*, wide=0pt]
  \item \textbf{Pre-allocation constraint:}
    In practice, there could be some pre-allocation of facilities and classrooms at some moment in some of the lecture rooms and laboratories.
    This model enforces that no course is assigned to such classrooms during the hours of unavailability.
    This is formulated as the following equalities
    \begin{align}
      x_{ijk} = 0 \quad \forall i \in I,\, \forall k \in K,\, \forall j \in J_{k}^{\text{unavail}}
      \eqfinp
    \end{align}

  \item \textbf{Infrastructure constraint:} 
    This first infrastructure constraint concerns with the capacity of a classroom.
    The following ensures that each student has a seat in a course:
    \begin{align}
      S_{i} x_{ijk} \leq C_{j} \quad \forall i \in I,\, \forall j \in J,\, \forall k \in K
      \eqfinp
    \end{align}
    We also require that a room is not shared between any two subjects at any moment:
    \begin{align}
      \sum_{i \in I} x_{ijk} \leq 1 \quad \forall j \in J,\, \forall k \in K
      \eqfinp
    \end{align}

    Beyond just capacity and unique utilization of a room, infrastructure matching uses metadata tags to completely segregate space allocation, explicitly forcing lecture courses to zero inside laboratory spaces and preventing laboratory courses from occupying standard lecture rooms.
    The constraints read as follow:
    \begin{align}
      x_{ijk} &= 0 \quad \forall i \in I_{\text{lec}},\, \forall j \in J_{\text{lab}},\, \forall k \in K 
      \eqfinv \\
      x_{ijk} &= 0 \quad \forall i \in I_{\text{lab}},\, \forall j \in J_{\text{lec}},\, \forall k \in K
      \eqfinp
    \end{align}

  \item \textbf{Unavailability constraint:} 
    Some of the instructors could be occupied with other commiments prior to the scheduling.
    The following equations force no teaching duties during the pre-occupied hours:
    \begin{align}
      x_{ijk} = 0 \quad \forall t \in T,\, \forall i \in I_{t},\, \forall k \in K_{t}^{\text{unavail}},\, \forall j \in J
      \eqfinp
    \end{align}
    The same is applied to the unavailable hours of the cohorts:
    \begin{align}
      x_{ijk} = 0 \quad \forall g \in G,\, \forall i \in I_{g},\, \forall j \in J,\, \forall k \in K_{g}^{\text{unavail}}
      \eqfinp
    \end{align}

  \item \textbf{Credit-hour fulfillment constraint:} 
    This constraint gurantees that each course has active hours equal to the requirement by the curriculum-design credits:
    \begin{align}\label{cons: credit-hour}
      \sum_{j \in J} \sum_{k \in K} x_{ijk} = H_{i} \quad \forall i \in I
      \eqfinp
    \end{align}

  \item \textbf{Fragmentation prevention constraint:} 
    We adopt the following equation to ensures that every class starts exactly once:
    \begin{align}
      \sum_{j \in J} \sum_{k \in K} y_{ijk} = 1 \quad \forall i \in I
      \eqfinp
    \end{align}
    We also need the following inequality to prevent courses from being broken into fragmented sub-sessions throughout the week:
    \begin{align}\label{cons: fragment forward}
      x_{i, j, k+\tau} \geq y_{ijk} \quad \forall i \in I,\, \forall j \in J,\, \forall k \in \{0, \dots, |K|-H_{i}-1\},\, \forall \tau \in \{0, \dots, H_{i}-1\}
      \eqfinp
    \end{align}

    Finally, since the time slots are indexed with integers, it is not immediate that two consecutive indices belong to the same day.
    The following constraint ensures that a course is contained within a single day:
    \begin{align}
      y_{ijk} = 0 \quad \forall i \in I,\, \forall j \in J,\, \forall k \in K^{\sharp}_{i}
      \eqfinv
    \end{align}
    where~$K^{\sharp}_{i} = \{ k \in K \mid D_{k} \neq D_{k+H_{i}-1} \ \text{or} \ k + H_{i} - 1 \geq |K|\}$.

  \item \textbf{Parallel session constraints:} 
    We consider in our model a special feature that a single instructor could conduct parallel sessions on designated courses.
    This is particularly practical in certain laboratory sessions where an instructor looks after concurrent groups along with some teaching assistants.

    First, we require that the parallel sessions are placed in different rooms by imposing the following inequalities:
    \begin{align}
      x_{s_1,j,k} + x_{s_2,j,k} \leq 1 \quad \forall t \in T,\, \forall j \in J,\, \forall k \in K,\, \forall I^{\text{par}}_{t} \in I^{\text{PAR}}_{t},\, \forall s_1, s_2 \in I^{\text{par}}_{t}
      \eqfinp
    \end{align}
    Moreover, to ensure that parallel sessions are actually parallel, we adopt the following constraints to synchronize their schedules:
    \begin{align}
      \sum_{j \in J} y_{s_1,j,k} = \sum_{j \in J} y_{s_2,j,k} \quad \forall t \in T,\, \forall k \in K,\, \forall I^{\text{par}}_{t} \in I^{\text{PAR}}_{t},\, \forall s_1, s_2 \in I^{\text{par}}_{t}
      \eqfinp
    \end{align}

    The subsequent inequalities under this constraint set are concerned with the assignment of an instructor.
    Outside of these special parallel arrangement, we must ensure that an instructor is either free or busy with a unique class:
    \begin{align}
      \sum_{i \in I_{t} \setminus I^{\text{par}}_{t}} \sum_{j \in J} x_{ijk} \leq 1 \quad \forall t \in T,\, \forall k \in K,\, \forall I^{\text{par}}_{t} \in I^{\text{PAR}}_{t}
      \eqfinp
    \end{align}
    The following final inequality warrants that no non-parallel course could be active along with the parallel sessions:
    \begin{align}
      \sum_{i \in I_{t}} \sum_{j \in J} x_{ijk} \leq \min_{I^{\text{par}}_{t} \in I^{\text{PAR}}_{t}}\abs{I^{\text{par}}_{t}}   \quad \forall t \in T,\, \forall k \in K
      \eqfinp
    \end{align}

  \item \textbf{Cohort conflict elimination constraint:} 
    This constraint eliminate class schedule conflicts for students within each cohort by restricting the total number of concurrent active courses, which reads
    \begin{align}
      \sum_{i \in I_{g}} \sum_{j \in J} x_{ijk} \leq 1 \quad \forall g \in G,\, \forall k \in K
      \eqfinp
    \end{align}

  \item \textbf{Lunch break constraint:} 
    This constraint rules out the possibility that a course overlaps the lunch break:
    \begin{align}
      y_{ijk} = 0 \quad \forall i \in I,\, \forall j \in J,\, \forall k \in K^{\flat}_{i}
      \eqfinv
    \end{align}
    where
    \begin{align}
      K^{\flat}_{i} = \Big\{ k \in K \;\Big|\; 0 \leq k \leq \abs{K} - H_{i},\: \{k, \dots, k + H_{i} - 1\} \cap K_{\text{lunch}} \neq \emptyset \Big\}
      \eqfinp
    \end{align}

  \item \textbf{Fatigue constriant:}
    To prevent fatigure, we require that a single instructor do not have two courses too close to gether to prevent fatigue:
    \begin{align*}
      \sum_{j \in J} ( y_{i,j,k} + y_{i',j,k'} ) \leq 1 \quad \forall t \in T, \forall i,i' \in I_{t}, i \neq i',\,  \forall k, k' \in K,\, k+H_{i} \leq k' \leq k+H_{i} + R - 1
      \eqfinp
    \end{align*}
    The same also applied to student cohorts:
    \begin{align*}
      \sum_{j \in J} ( y_{i,j,k} + y_{i',j,k'} ) \leq 1 \quad \forall g \in G,\, \forall i,i' \in I_{g}, i \neq i',\,  \forall k, k' \in K, k+H_{i} \leq k' \leq k+H_{i} + R - 1
      \eqfinp
    \end{align*}
\end{itemize}

\section{Implementaion: A toy example}
\label{sec: toy example}

In this section, we demonstrate an implementation of our model with a small-scale toy example using {\ttfamily Gurobi Optimizer} via {\ttfamily Python} interface {\ttfamily gurobipy}.
The problem instance has 11 subjects, of which 8 are lecture subjects and the remaining 3 are laboratory subjects.
There are 3 instructors, 4 student cohorts, 3 classrooms, and no parallel classes.
The timeslots available for assignment consists of six periods per day from Monday to Friday, each of which is one-hour long.
Every subjects must be scheduled according to all the rules presented in~\S\ref{subsec: constraints}.
The resulting schedule illustrates not only feasibility with respect to institutional resources, but also a more realistic timetable structure for instructors and students.

\subsection{Problem profile}

We start by presenting the full details of this toy example to illustrate as much as possible the construction and implementation of our BLP model.
The dimension of our sample data is summarized in Table~\ref{tab: toy-basic-size}.
The characteristics of each cohorts and classrooms are also depicted in Tables~\ref{tab: toy cohorts} and~\ref{tab: toy classrooms}, respectively.

\begin{table}[H]
  \centering
  \begin{minipage}[t]{.5\textwidth}
    \caption{Problem instance profile}
    \label{tab: toy-basic-size}
    \begin{tblr}{
        colspec = {l r},
        row{1}={font=\bfseries},
        hline{1,Z}={0.8pt},
        hline{2}={0.5pt},
        rowsep=3pt,
        colsep=6pt
      }
      Item & Dimension \\
      Subjects ($I$) & 11 \\
      \quad - Lecture subjects ($I_{\text{lec}}$) & 8 \\
      \quad - Laboratory subjects ($I_{\text{lab}}$) & 3 \\
      Instructors ($T$) & 3 \\
      Student cohorts ($G$) & 4 \\
      Temporal slots ($K$) & 30 \\
      Classrooms ($J$) & 3 \\
      \quad - Lecture rooms ($J_{\text{lec}}$) & 2 \\
      \quad - Laboratory ($J_{\text{lab}}$) & 1 \\
      Parallel subject pairs & 0 \\
    \end{tblr}
  \end{minipage}
  \begin{minipage}[t]{.45\textwidth}
    \centering
    \caption{Cohort data}
    \label{tab: toy cohorts}
    \begin{tblr}{
        colspec={cc},
        row{1}={font=\bfseries},
        hline{1,Z}={0.8pt},
        hline{2}={0.5pt},
        rowsep=3pt,
        colsep=6pt
      }
      Cohort & Number of students \\
      C1 & 24 \\
      C2 & 22 \\
      C3 & 26 \\
      C4 & 22 \\
    \end{tblr}
    \vskip 1.2em
    \caption{Classroom data}
    \label{tab: toy classrooms}
    \begin{tblr}{
        colspec={ccc},
        row{1}={font=\bfseries},
        hline{1,Z}={0.8pt},
        hline{2}={0.5pt},
        rowsep=3pt,
        colsep=6pt
      }
      Room & Type & Capacity \\
      LEC-A & Lecture & 110 \\
      LEC-B & Lecture & 60 \\
      LAB-A & Lab & 55 \\
    \end{tblr}
  \end{minipage}
\end{table}

Next we present the details regarding the subjects, class types, class length, the associated instructors, the associated student cohorts, and the enrollment sizes used in this example.
This complete information is shown below in the Table~\ref{tab: toy subject profiles}.
\begin{longtblr}[
  caption={Subject and Cohort profiles},
  label={tab: toy subject profiles},
  ]{
    colspec={c X[l] l c c c c c c c},
    row{1,2}={font=\bfseries},
    hline{1,Z}={0.8pt},
    hline{2}={0.5pt}, 
    hline{3}={0.5pt}, 
    rowsep=3pt,
    colsep=6pt
  }
  \SetCell[r=2]{c} ID & \SetCell[r=2]{l} Subject & \SetCell[r=2]{l} Type & \SetCell[r=2]{c} Hours & \SetCell[r=2]{c} Instructor & \SetCell[c=4]{c} Cohorts & & & & \SetCell[r=2]{r} Enrolled \\
                      &            &            &       &         & C1 & C2 & C3 & C4 &         \\
  S01 & Calculus I & Lecture    & 2     & T1      & $\checkmark$ & $\checkmark$ & & & 46 \\
  S02 & Linear Algebra & Lecture & 2    & T1      & $\checkmark$ & & $\checkmark$ & & 50 \\
  S03 & Discrete Mathematics & Lecture & 2 & T2   & & $\checkmark$ & & $\checkmark$ & 44 \\
  S04 & Optimization & Lecture & 2   & T3      & & & $\checkmark$ & $\checkmark$ & 48 \\
  S05 & Probability & Lecture  & 2     & T2      & $\checkmark$ & & & $\checkmark$ & 46 \\
  S06 & Statistics & Lecture   & 2     & T2      & & $\checkmark$ & $\checkmark$ & & 48 \\
  S07 & Scientific Python & Lecture & 3 & T1   & $\checkmark$ & $\checkmark$ & $\checkmark$ & & 72 \\
  S08 & Mathematical Modeling & Lecture & 2 & T3 & $\checkmark$ & $\checkmark$ & $\checkmark$ & $\checkmark$ & 94 \\
  S09 & Python Lab & Laboratory & 2  & T1      & $\checkmark$ & $\checkmark$ & & & 46 \\
  S10 & Statistics Lab & Laboratory & 2 & T2   & & $\checkmark$ & $\checkmark$ & & 48 \\
  S11 & Optimization Lab & Laboratory & 2 & T3   & & & $\checkmark$ & $\checkmark$ & 48 \\
\end{longtblr}

Finally, we assume some unavailability restrictions on instructors, classrooms and student cohorts.
These constraints are described in the following Table~\ref{tab: toy unavailability}.

\begin{table}[H]
  \centering
  \caption{Unavailability restrictions}
  \label{tab: toy unavailability}
  \begin{tblr}{
      colspec={c c l},
      row{1}={font=\bfseries},
      hline{1,Z}={0.8pt},
      hline{2}={0.5pt},
      rowsep=3pt,
      colsep=6pt
    }
    Entity & Unavailable period & Reason \\
    Teacher T1 & Fri-P5, Fri-P6 & Teacher unavailable \\
    Teacher T2 & Mon-P1 & Teacher unavailable \\
    Teacher T3 & Wed-P6 & Teacher unavailable \\
    Room LAB-A & Tue-P6 & Lab maintenance \\
    Room LEC-B & Thu-P2 & Room reserved \\
    Cohort C4 & Fri-P6 & Cohort activity \\
  \end{tblr}
\end{table}

\subsection{Optimization results}

In this subsection, we outline the solution of our BLP model obtained through {\ttfamily Gurobi Optimizer} as explained at the beginning of~\S\ref{sec: toy example}.
The optimal timetable contains 11 scheduled class blocks and 23 hours of classroom utilization in total.  
The model produces no room-type, teacher, cohort, or lunch-period violation.  
We summarize the optimal solution in Table~\ref{tab: toy solution-summary}.
\begin{table}[H]
  \centering
  \caption{Solution summary}
  \label{tab: toy solution-summary}
  \begin{tblr}{
      colspec={l r},
      row{1}={font=\bfseries},
      hline{1,Z}={0.8pt},
      hline{2}={0.5pt},
      rowsep=3pt,
      colsep=6pt
    }
    Item & Value \\
    Objective value & 2915.206 \\
    Scheduled class blocks & 11 \\
    Classroom utilization hours & 23 \\
    Seat efficiency & 78.25\% \\
  \end{tblr}
\end{table}

Extracted from the optimal solution of our model, the variable values are piped into our own script to make a readable timetables for each instructors.
These are shown in Tables~\ref{tab: toy table T1}--\ref{tab: toy table T3}.

\begin{table}[H]
  \centering
  \caption{Schedule for Instructor T1}
  \label{tab: toy table T1}
  \begin{tblr}{
      colspec={c c l c l},
      row{1}={font=\bfseries},
      hline{1,Z}={0.8pt},
      hline{2}={0.5pt},
      rowsep=3pt,
      colsep=6pt
    }
    Day & Time & Subject & Room & Cohorts \\
    Tue & 10:30--12:30 & S09 Python Lab & LAB-A & C1, C2 \\
    Thu & 08:30--11:30 & S07 Scientific Python & LEC-A & C1, C2, C3 \\
    Thu & 13:30--15:30 & S02 Linear Algebra & LEC-B & C1, C3 \\
    Fri & 10:30--12:30 & S01 Calculus I & LEC-B & C1, C2 \\
  \end{tblr}
\end{table}

\begin{table}[H]
  \centering
  \caption{Schedule for Instructor T2}
  \label{tab: toy table T2}
  \begin{tblr}{
      colspec={c c l c l},
      row{1}={font=\bfseries},
      hline{1,Z}={0.8pt},
      hline{2}={0.5pt},
      rowsep=3pt,
      colsep=6pt
    }
    Day & Time & Subject & Room & Cohorts \\
    Mon & 13:30--15:30 & S10 Statistics Lab & LAB-A & C2, C3 \\
    Wed & 10:30--12:30 & S03 Discrete Mathematics & LEC-B & C2, C4 \\
    Wed & 13:30--15:30 & S05 Probability & LEC-B & C1, C4 \\
    Fri & 13:30--15:30 & S06 Statistics & LEC-B & C2, C3 \\
  \end{tblr}
\end{table}

\begin{table}[H]
  \centering
  \caption{Schedule for Instructor T3}
  \label{tab: toy table T3}
  \begin{tblr}{
      colspec={c c l c l},
      row{1}={font=\bfseries},
      hline{1,Z}={0.8pt},
      hline{2}={0.5pt},
      rowsep=3pt,
      colsep=6pt
    }
    Day & Time & Subject & Room & Cohorts \\
    Tue & 08:30--10:30 & S11 Optimization Lab & LAB-A & C3, C4 \\
    Tue & 13:30--15:30 & S08 Mathematical Modeling & LEC-A & C1, C2, C3, C4 \\
    Fri & 08:30--10:30 & S04 Optimization & LEC-B & C3, C4 \\
  \end{tblr}
\end{table}

To sum up the obtained solution, the overall optimal schedule is collected in the Table~\ref{tab: toy opt sched} below.
\begin{longtblr}[
  caption={Optimal schedule ordered by the time sequence},
  label={tab: toy opt sched}
  ]{
    colspec={c X[l] c c X[l] c c c},
    row{1}={font=\bfseries},
    hline{1,Z}={0.8pt},
    hline{2}={0.5pt},
    rowsep=3pt,
    colsep=6pt
  }
  ID & Subject & Teacher & Room & Cohorts & Day & Time & Empty seats \\
  S10 & Statistics Lab & T2 & LAB-A & C2, C3 & Mon & 13:30--15:30 & 7 \\
  S11 & Optimization Lab & T3 & LAB-A & C3, C4 & Tue & 08:30--10:30 & 7 \\
  S09 & Python Lab & T1 & LAB-A & C1, C2 & Tue & 10:30--12:30 & 9 \\
  S08 & Mathematical Modeling & T3 & LEC-A & C1, C2, C3, C4 & Tue & 13:30--15:30 & 16 \\
  S03 & Discrete Mathematics & T2 & LEC-B & C2, C4 & Wed & 10:30--12:30 & 16 \\
  S05 & Probability & T2 & LEC-B & C1, C4 & Wed & 13:30--15:30 & 14 \\
  S07 & Scientific Python & T1 & LEC-A & C1, C2, C3 & Thu & 08:30--11:30 & 38 \\
  S02 & Linear Algebra & T1 & LEC-B & C1, C3 & Thu & 13:30--15:30 & 10 \\
  S04 & Optimization & T3 & LEC-B & C3, C4 & Fri & 08:30--10:30 & 12 \\
  S01 & Calculus I & T1 & LEC-B & C1, C2 & Fri & 10:30--12:30 & 14 \\
  S06 & Statistics & T2 & LEC-B & C2, C3 & Fri & 13:30--15:30 & 12 \\
\end{longtblr}

The computational result shows that the proposed formulation is able to produce a feasible and interpretable timetable for the toy instance. 
All subjects are assigned to compatible rooms and valid time blocks, while the teacher, cohort, room, availability, lunch-break, and rest-period requirements are simultaneously satisfied. 
This confirms that the model can incorporate practical scheduling considerations without changing the main assignment structure, and provides a useful basis for extending the formulation to larger and more realistic class-scheduling problems for the next section.

\section{Implementation: A case study at KMUTT}
\label{sec: case study}

In this section, we apply our proposed model from~\S\ref{subsec: model} to a case study of Department of Mathematics, KMUTT, for its scheduling in the Semester~1/2025.
As always, the binary linear programming (BLP) framework was solved with the~{\ttfamily Gurobi Optimizer}.
The implementation was written in {\ttfamily Python} through the {\ttfamily gurobipy} interface.
The experiment were processed on a standard laptop with Intel Core Ultra 5 135U with 16GB DDR5 RAM.

\subsection{Problem profile}

The complexity and combinatorial scale of our model are directly related to the base dimensions of the institution's profile.
The characteristics of the dataset are formally summarized in Table \ref{tab:dataset_profile}.
\begin{table}[H]
\centering
\caption{Specifications of the evaluated dataset instance}
\label{tab:dataset_profile}
\begin{tblr}{
  colspec = {l r},
  hline{1,Z} = {0.08em}, 
  hline{2} = {0.05em},   
}
\textbf{Component} & \textbf{Dimension} \\
Subjects ($I$) & 55 \\
Instructors ($T$) & 29 \\
Student cohorts ($G$) & 126 \\
Sum total of unique enrolled students & 489 \\
Temporal slots ($K$) & 40 \\
Classrooms ($J$) & 55 \\
\quad -- Lecture Rooms & 52 \\
\quad -- Laboratory Rooms & 3 \\
\end{tblr}
\end{table}

\subsection{Optimization results}

After the data pre-processing, we obtain the needful sets and parameters for the BLP model.
This process yields a high-dimensional combinatorial problem matrix. 
The internal complexity of the model is reflected by the pre-processed constraints rows count and active binary column structures generated prior to the solver convergence.
Table~\ref{tab:model_performance} logs the computational parameters recorded during execution.
\begin{table}[H]
\centering
\caption{The solver log}
\label{tab:model_performance}
\begin{tblr}{
  colspec = {l r},
  hline{1,Z} = {0.08em},
  hline{2} = {0.05em},
}
\textbf{Execution parameter} & \textbf{Value} \\
Solver termination status & Optimal (status code: 2) \\
Minimum objective value & 19,648.00 \\
Original constraints rows & 621,590 \\
Original binary decision variables & 242,000 \\
Presolved / reduced constraints rows & 41,308 \\
Presolved / reduced binary decision variables & 62,380 \\
Cpu time & 80.36180 seconds \\
\end{tblr}
\end{table}

Moreover, the  Table~\ref{tab:model_performance} indicates that the global optimality is reached by the solver in approximately 80.36 seconds. 
Given the profound scale of the underlying matrix, originally incorporating more than 621,000 active rows of structural constraints and over 242,000 active variables, the computation time remains operationally tractable for institutional scheduling.

\subsection{Resource utilization and comparison with manual schedule}

The primary objective function drives the optimization toward an optimal allocation by minimizing a squared penalty. 
This design selectively penalizes the slacks from large classrooms, influencing the solution to avoid seating small student cohorts in oversized lecture rooms. 
However, this squared objective may distort the direct understanding of the number of vacant seats itself.
To highlight the efficiency of the proposed method, we also compute the total actual empty seats from our optimal schedule, which equals to~934 seats across horizon of schuduling.

We also compare our optimal schedule obtained from our proposed model with the schedule derived manually by hand from the Secretary of Department of Mathematics, KMUTT, who has been responsible for this task for over five years.
It turns out that this manually constructed schedule exhibits as high as 1,525 empty seats across horizon.
This means we are able to reduce~38.75\% of empty seats with our proposed BLP mode under this case study.

Below, the Table \ref{tab:utilization_summary} reports several aspects of allocation metrics including the total actual unused seats.

\begin{table}[htbp]
\centering
\caption{Resource allocation and comparison with the manual schedule}
\label{tab:utilization_summary}
\begin{tblr}{
  colspec = {l r},
  hline{1,Z} = {0.08em},
  hline{2} = {0.05em},
}
\textbf{Performance evaluation metric} & \textbf{Value} \\
Averaged seat occupation rate (Proposed model) & 63.21\% \\
Total actual empty seats across horizon (Proposed model)  & 934 seats \\
Total actual empty seats across horizon (Manual method)  & 1,525 seats \\
\textbf{Average empty seats per scheduled hour (Proposed model)} & \textbf{16.98 seats / hour} \\
\textbf{Improvement over the Manual method} & \textbf{38.75\%}
\end{tblr}
\end{table}

The optimal solution establishes that, on average, a classroom is~63.21\% full.
This percentage attributes to the presence of many small specialized subjects with few students, which must still be assigned to lecture rooms that are substantially larger.

When analyzing the linear room metrics, the network maintains an average of only 16.98 empty seats during any moment of active hours. 
The optimal solution in this case study proves that the quadratic formulation successfully restricted massive space slacks. 

\small 
\section*{Acknowledgments}

The authors are grateful to Angkool Wangwongchai, the secretary of the Department of Mathematics, King Mongkut's University of Technology Thonburi, for providing the experience, insight and all the data used in the study.

\section*{Author contributions}

\begin{itemize}[label=$\circ$, leftmargin=*]
  \item Supervision: Parin Chaipunya; 
  \item Methodology: Parin Chaipunya;
  \item Formal analysis: Ohmchana Klinsod;
  \item Investigation: Ohmchana Klinsod, Parin Chaipunya;
  \item Data curation: Parin Chaipunya;
  \item Data preparation and processing: Ohmchana Klinsod;
  \item Software -- Ohmchana Klinsod;
  \item Writing -- Original draft: Ohmchana Klinsod, Parin Chaipunya; 
  \item Writing -- Review \& editing, Ohmchana Klinsod, Parin Chaipunya;
\end{itemize}

\section*{Data availability}

The data underlying this study were provided by the authors' institution under internal use and are not publicly available.

\bibliography{class_sched.bib}

\end{document}